\documentstyle[amsfonts,12pt]{article}

\newtheorem{thm}{Teorem}[section]

\newtheorem{lem}[thm]{\quad Lemma}

\newtheorem{defi}[thm]{\quad Definition}
\newtheorem{teo}[thm]{\quad Theorem}

\newtheorem{ex}[thm]{\quad Example}

\begin{document}

\title{A Unified Theory on Some Basic
Topological Concepts}
\author{T. Hatice Yalva\c{c}}
\maketitle

\begin{abstract}
Several mathematicians, including myself, have studied some unifications in
general topological spaces as well as in fuzzy topological spaces. For
instance in our earlier works, using operations on topological spaces, we
have tried to unify some concepts similar to continuity, openness,
closedness of functions, compactness, filter convergence, closedness of
graphs, countable compactness and Lindel\"{o}f property. In this article, to
obtain further unifications, we will study $\varphi_{1,2}$-compactness and
relations between $\varphi_{1,2}$-compactness, filters and $\varphi_{1,2}$%
-closure operator.
\end{abstract}

\section{Introduction}

Several unifications have been studied in [1, 9, 10, 11, 12, 13, 14, 15, 16,
17, 18, 19, 20, 21, 22]. Some unifications were studied in [10] and [12] by
using operations for fuzzy topological spaces. It was claimed there and it
can be easly seen that most of the definitions and results can be applied to
topological spaces. As far as possible we do not repeat definitions related
of known concepts. Because one aim of us to reduce the confusions caused by
so much definitions. However, many such definitions will be clear from the
special operations considered.

In a topological space $(X, \tau)$ int, cl, scl etc. will stand for
interior, closure, semi-closure operations so on and $A^{o}$, $\overline{A}$
will stand for the interior of A, the closure of A for a subset A of X
respectively.\newline

\begin{defi}
Let $(X, \tau)$ be a topological space. A mapping $\varphi : P(X)
\rightarrow P(X)$ is called an operation on $(X, \tau)$ if $A^{o} \subset
\varphi (A)$ for all $A \in P(X)$ and $\varphi (\emptyset) = \emptyset$.
\end{defi}

The class of all operations on a topological space $(X, \tau)$ will be
denoted by $O(X, \tau)$.

A partial order "$\leq$" on $O(X, \tau)$ is defined as $\varphi_{1} \leq
\varphi_{2} \Leftrightarrow \varphi_{1}(A) \subset \varphi_{2}(A)$ for each $%
A \in P(X)$.

An operation $\varphi \in O(X, \tau)$ is called monotonous if $\varphi (A)
\subset \varphi (B)$ whenever $A \subset B $ $(A, B \in P(X))$.\newline

\begin{defi}
Let $\varphi \in O(X, \tau)$ and $A,B\subset X$. $A$ is called $\varphi$%
-open if $A\subset \varphi (A)$. $B$ is called $\varphi$-closed if $%
X\setminus B$ is $\varphi$-open. Operation $\tilde \varphi \in O(X, \tau)$
is called the dual operation of $\varphi$ if $\tilde \varphi (A)= X\setminus
\varphi (X\setminus A)$ for each $A \in P(X)$.
\end{defi}

If $\varphi$ is monotonous, then the family of all $\varphi$-open sets is a
supratopology (${\cal U} \subset P(X)$ is a supratopology on X means that $%
\emptyset \in {\cal U}$, $X\in {\cal U}$ and ${\cal U}$ is closed under
arbitrary union [2]).

Let $(X,\tau)$ be a topological space, $\varphi \in O(X,\tau)$, ${\cal U}
\subset P(X), x \in X$. We use the following notations.\newline
\[
{\cal U} (x) =\{U: x\in U\in {\cal U} \},
\]
\[
\qquad \quad\quad\quad\varphi O(X)=\{U: U \subset X, U \quad is \quad
\varphi-open\},
\]
\[
\qquad \quad\quad\quad\varphi C(X)=\{K: K \subset X, K \quad is \quad
\varphi-closed\},
\]
\[
\varphi O(X,x)=\{U: U \in \varphi O(X), x \in U \}.
\]
\[
{\cal N}({\cal U},x)=\{N: N \subset X \quad and \quad there \quad exists
\quad a \quad U \in {\cal U} (x) \quad such \quad that \quad U \subset N\}.
\]

\begin{defi}
Let $\varphi \in O(X, \tau)$, ${\cal U} \subset P(X)$. $\varphi$ is called
regular with respect to (shortly w.r.t.) ${\cal U}$ if $x \in X$ and $U,V
\in {\cal U} (x)$, there exists an $W\in {\cal U} (x)$ such that $\varphi
(W) \subset \varphi (U) \cap \varphi (V)$.
\end{defi}

For any operation $\varphi\in O(X, \tau)$, $\tau\subset \varphi O(X)$, and $%
X, \emptyset$ are both $\varphi$-open and $\varphi$-closed.\newline

\begin{defi}
Let $\varphi_{1}, \varphi_{2} \in O(X,\tau) , A \subset X $.

a) $x \in \varphi_{1,2} intA \Leftrightarrow there \quad exists \quad an
\quad U \in \varphi_{1} O(X,x) \quad such \quad that \quad$ \newline
$\varphi_{2}(U) \subset A $.

b) $x \in \varphi_{1,2} clA \Leftrightarrow for\quad each \quad U\in
\varphi_{1} O(X,x), \quad \varphi_{2} (U) \cap A \neq \emptyset$.

c) A is $\varphi_{1,2}$-$open \Leftrightarrow A \subset \varphi_{1,2} intA $.

d) A is $\varphi_{1,2}$-$closed \Leftrightarrow \varphi_{1,2} clA \subset A $%
.
\end{defi}

If $A\subset B$ then $\varphi_{1,2} intA \subset \varphi_{1,2} intB$.
Clearly for any set $A, X \setminus \varphi_{1,2} intA = \varphi_{1,2} cl(X
\setminus A)$ and  A is $\varphi_{1,2}$-open iff $X \setminus A \quad
is\quad \varphi_{1,2}$-$closed $.

$\varphi_{1,2} O(X)$\quad $(\varphi_{1,2} C(X))$ will stand for the family
of all $\varphi_{1,2}$-open subsets ($\varphi_{1,2}$-closed subsets) of X.

\begin{teo}
([12]). Let $\varphi_{1} , \varphi_{2} \in O(X, \tau) $.

a) $\varphi_{1,2} O(X) $ is a supratopology on X.

b) If $\varphi_{2}$ is regular w.r.t. $\varphi_{1} O(X)$ then $\varphi_{1,2}
O(X)$ is a topology on X and a subset K of X is closed w.r.t. this topology
iff $\varphi_{1,2} clK \subset K$.

c) If $\varphi_{2}$ is regular w.r.t. $\varphi_{1} O(X)$ and ($\varphi_{2}
\geq \imath \quad or \quad \varphi_{2} \geq \varphi_{1}$) then $%
\varphi_{1,2} O(X)$ is a topology on X and a set K is closed w.r.t. this
topology iff $\varphi_{1,2} cl K =K$ (here $\imath $ is the idendity
operation).
\end{teo}

\begin{ex}
Let the following operations be defined on a topological space $(X, \tau)$.

$\varphi_{1} =int$, \quad $\varphi_{2} =cloint$, \quad $\varphi_3 =cl$,
\quad $\varphi_4 =scl$, \quad $\varphi_5 =\imath$, \quad $\varphi_6 =intocl$.%
\newline
$\varphi_{1} \leq \varphi_{2} \leq \varphi_3$, \quad $\varphi_{1} \leq
\varphi _5 \leq \varphi _4 \leq \varphi _3$, \quad $\varphi_{1} \leq \varphi
_6 \leq \varphi _4$.\newline
$\varphi_{1} O(X)=\tau$, \quad $\varphi_{2} O(X)=SO(X)=$ the family of
semi-open sets.\newline
$\varphi_3 O(X)=\varphi_5 O(X)=\varphi_4 O(X)=P(X)=$ power set of X.\newline
$\varphi_6 O(X)=PO(X)=$ the family of pre-open sets.\newline
$\varphi_{1,3}O(X)=\tau _{\theta}=$ the topology of all $\theta$-open sets.%
\newline
$\varphi_{2,4}O(X)=S\theta O(X)=$ the family of semi-$\theta$-open sets.%
\newline
$\varphi_{1,6}O(X)=\tau_s =$ the semi regularization topology of X. It is,
the topology with the base RO(X) which consists of regular open sets = the
family of $\delta$-open sets.\newline
$\varphi_{2,3}O(X)= \theta SO(X)=$ the family of all $\theta$-semi-open sets.%
\newline
$\varphi_{1}, \varphi_3$ ($\varphi_{2}, \varphi_6$) are the dual operations
of each other. $\varphi_{2}$, $\varphi_3$, $\varphi _4$, $\varphi_5$, $%
\varphi_6$ are regular w.r.t. $\varphi_{1} O(X)$.

SC(X) (PC(X), RC(X), S$\theta$C(X), $\theta$SC(X) respectively) will stand
for the family semi-closed (pre-closed, regular closed, semi-$\theta$%
-closed, $\theta$-semi-closed) sets. \newline
SR(X)=SO(X)$\cap$SC(X)= the family of semi-regular sets.
\end{ex}

For operations $\varphi_{1}$, $\varphi_{2} \in O(X, \tau)$, clearly if $%
\varphi_{1} $ is monotonous and $\varphi_{2} =\imath$ for $\varphi_{1}
,\varphi_{2} \in O(X, \tau)$ then $\varphi_{1,2} O(X)=\varphi_{1} O(X)$ and $%
\varphi_{1,2} C(X)=\varphi_{1} C(X)$. If $\varphi_{1} O(X)$ is a topology
and $\varphi_{2}$ is monotonous then $\varphi_{2}$ is regular w.r.t. $%
\varphi_{1} O(X)$, so $\varphi_{1,2} O(X)$ is always a topology.

\section{$\protect\varphi _{1,2}$-closure Operator and Filters}

Along of the paper it will be accepted that operations $\varphi _i$, $%
i=1,2,...$ are defined on a topological space $(X, \tau)$.

\begin{lem}
If ($\varphi_{2} \geq \varphi_{1}$ or $\varphi_{2} \geq \imath$) then $%
\varphi_{1} O(X)\subset \varphi_{2} O(X)$. But the converse is not true.
\end{lem}

\begin{ex}
Let $\varphi_{1} =cloint$, $\varphi_{2} =$semi-int be defined on ${\Bbb R}$
with the usual topology. For $A=(0,1]$, $\varphi_{1} (A)=[0,1]$, $%
\varphi_{2} (A)=(0,1]$ and $\varphi_{1} (A) \not\subset \varphi_{2} (A)$.
i.e. $\varphi_{1} \not\leq \varphi_{2}$.

For the set of rational numbers ${\Bbb Q}$, ${\Bbb Q}\not\subset \varphi_{2}
({\Bbb Q})=\emptyset$. i.e. $\varphi_{2} \not\geq \imath$. But $\varphi_{1}
O({\Bbb R})=SO({\Bbb R})=\varphi_{2} O({\Bbb R})$.
\end{ex}

\begin{teo}
Let ${\cal B} =\{\varphi_{2} (U):U\in \varphi_{1} O(X)\}$.

a) If $\varphi_{2}$ is regular w.r.t. $\varphi_{1} O(X)$ and $\varphi_{1}
O(X) \subset \varphi_{2} O(X)$ then $\varphi_{1,2}$-cl operator defines the
same (pre)topology given in Theorem 1.5 (c). Let $\tau _{\varphi_{1,2}} $ be
stand for this topology. $A\subset \varphi_{1,2} cl A \subset \tau
_{\varphi_{1,2}} clA$ for any subset A of X.\newline

b) If $\varphi_{2}$ is regular w.r.t. $\varphi_{1} O(X)$, $\varphi_{1} O(X)
\subset \varphi_{2} O(X)$ and ${\cal B} \subset \varphi_{1,2} O(X)$, then $%
\varphi_{1,2}$-cl operator is a Kuratowski closure operator and $%
\varphi_{1,2} cl A= \tau _{\varphi_{1,2}} clA$ for any subset A of X.\newline
\end{teo}

\begin{defi}
([20]). Let ${\cal F}$ be a filter (or filterbase) in $(X, \tau)$ and $a\in X
$. ${\cal F}$ is said to be:\newline

a) $\varphi_{1,2}$-accumulates to a if $a\in \cap \{\varphi_{1,2} cl F :
F\in {\cal F}\}$.\newline

b) $\varphi_{1,2}$-converges to a if for each $U\in \varphi_{1} O(X,a)$,
there exists an $F\in {\cal F}$ such that $F\subset \varphi_{2} (U)$.
\end{defi}

\begin{teo}
([20]). 1) A filterbase ${\cal F} _b$ $\varphi_{1,2}$-accumulates ($%
\varphi_{1,2}$-converges) to a iff filter generated by ${\cal F} _b$ $%
\varphi_{1,2}$-accumulates ($\varphi_{1,2}$-converges) to a.\newline

2) If $\varphi_{2}$ is monotonous, we can get the family ${\cal N}
(\varphi_{1} O(X),a)$ instead of $\varphi_{1} O(X,a)$ in the above
definitions.\newline

3) A filter ${\cal F}$ $\varphi_{1,2}$-converges to a iff $\{\varphi_{2}
(U): U\in \varphi_{1} O(X,a)\}\subset {\cal F}$.\newline

4) If ${\cal F}$ $\varphi_{1,2}$-converges to a then ${\cal F}$ $%
\varphi_{1,2}$-accumulates to a.\newline

5) Let ${\cal F} \subset {\cal F} ^{^{\prime}}$ for the filters ${\cal F} $
and ${\cal F} ^{^{\prime}}$

a) If ${\cal F} ^{^{\prime}}$ $\varphi_{1,2}$-accumulates to a, then ${\cal F%
} $ $\varphi_{1,2}$-accumulates to a.

b) If ${\cal F} $ $\varphi_{1,2}$-converges to a then ${\cal F} ^{^{\prime}}$
$\varphi_{1,2}$-converges to a.\newline

6) If $\varphi_{2}$ is regular w.r.t. $\varphi_{1} O(X)$ then a filter $%
{\cal F}$ $\varphi_{1,2}$-accumulates to $a\in X$ iff there exists a filter $%
{\cal F} ^{^{\prime}}$ such that ${\cal F} \subset {\cal F} ^{^{\prime}}$
and ${\cal F} ^{^{\prime}}$ $\varphi_{1,2}$-converges to a.\newline

7) If ${\cal F}$ is a maximal filter which $\varphi_{1,2}$-accumulates to $%
a\in X$, then ${\cal F}$ $\varphi_{1,2}$-converges to a.\newline

8) If $\varphi_{1}^{^{\prime}} O(X) \subset \varphi_{1} O(X)$ and $%
\varphi_{2}^{^{\prime}} \geq \varphi_{2}$ for the operations $\varphi_{1}$, $%
\varphi_{2}$, $\varphi_{1}^{^{\prime}}$, $\varphi_{2}^{^{\prime}} \in O(X,
\tau)$, then a filter (or a filterbase) ${\cal F}$ $\varphi_{1,2}^{^{\prime}}
$-accumulates ($\varphi_{1,2}^{^{\prime}}$-converges) to a whenever ${\cal F}
$ $\varphi_{1,2}$-accumulates ($\varphi_{1,2}$-converges) to a.
\end{teo}

{\bf Proof.} 6) Let $\varphi_{2}$ be regular w.r.t. $\varphi_{1} O(X)$ and $%
{\cal F}$ $\varphi_{1,2}$-accumulates  to a. Then the family ${\cal F} _{b}=
\{\varphi_{2} (U)\cap F: U\in \varphi_{1} O(X,a),F\in {\cal F}\}$ is a
filterbase and the filter ${\cal F} ^{^{\prime}}$ generated by ${\cal F} _b$
is finer than ${\cal F}$ and it $\varphi_{1,2}$-converges to a. The other
part of the proof is clear.\newline

7) Let ${\cal F}$ be a maximal filter and $\varphi_{1,2}$-accumulates to $%
a\in X$. For each $U\in \varphi_{1} O(X,a)$, ${\cal F}_{b_u}=\{\varphi_{2}
(U)\cap F: F\in {\cal F}\}$ is a filterbase. The filter ${\cal F} _{u}$
generated by ${\cal F} _{b_u}$ is finer than ${\cal F}$. So, for each $U\in
\varphi_{1} O(X,a)$, ${\cal F} = {\cal F}_u$. Now it is clear that ${\cal F}$
$\varphi_{1,2}$-converges to a.

A space $(X, \tau)$ is called $\varphi_{1,2}$-$T_2$ if for each $x,y\in
X(x\neq y)$ there are $\varphi_{1}$-open sets $U_x$ and $U_y$ such that $%
x\in U_x$, $y\in U_y$ and $\varphi_{2} (U_x) \cap \varphi_{2} (U_y)
=\emptyset$ ([12],[20]).

\begin{teo}
([20]) Let $(X, \tau)$ be $\varphi_{1,2}$-$T_2$ space. If a filter ${\cal F}$
$\varphi_{1,2}$-converges to some point $a\in X$ and $\varphi_{1,2}$%
-accumulates to some point $b\in X$ then $a=b$.
\end{teo}

{\bf Proof.} Let's accept that a filter ${\cal F}$ be $\varphi_{1,2}$%
-convergent to $a$ and $\varphi_{1,2}$-accumulate to $b$ and $a\neq b$. Then
there exists $U\in \varphi_{1} O(X,a)$ and $V\in \varphi_{1} O(X,b)$ such
that $\varphi_{2} (U) \cap \varphi_{2} (V) =\emptyset$. But $\varphi_{2}
(U)\in {\cal F}$ and $F\cap \varphi_{2} (V)\neq \emptyset$ for each $F\in
{\cal F}$. It must be $\varphi_{2} (U) \cap \varphi_{2} (V)\neq \emptyset$.
This contradiction completes the proof.

\begin{ex}
Let $a\in X$ and ${\cal F}$ be a filter in $(X, \tau)$.\newline

a) Let $\varphi_{1} = cloint$, \quad $\varphi_{2} = cl$.

${\cal F}$ $\varphi_{1,2}$-converges ($\varphi_{1,2}$-accumulates) to a iff $%
{\cal F}$ rc-converges (rc-accumulates) to a since $\{\overline V : V\in
\tau, x\in \overline V\}=\{\overline U : x\in U\in SO(X)\}$ [8] iff ${\cal F}
$ s-converges (s-accumulates) to a [5]. \newline

b) Let $\varphi_{1} = int$, \quad $\varphi_{2} =cl$.

${\cal F}$ $\varphi_{1,2}$-converges ($\varphi_{1,2}$-accumulates) to a iff $%
{\cal F}$ r-converges (r-accumulates) to a [7] iff ${\cal F}$ almost
converges to a (a is an almost adherent point of ${\cal F}$ )[4].

$(X, \tau)$ is $\varphi_{1,2}$-$T_2$ iff $(X, \tau)$ is Urysohn.\newline

c) Let $\varphi_{1} =int$, \quad $\varphi_{2} =\imath$.

${\cal F}$ $\varphi_{1,2}$-converges to a iff ${\cal F}$ converges to a in $%
(X, \tau)$.

$(X, \tau)$ is $\varphi_{1,2}$-$T_2$ iff $(X, \tau)$ is Hausdorff.

It is well known that in Hausdorff spaces, a convergent filter can not have
more than one accumulation point ([6], page 220).
\end{ex}

\begin{teo}
([20])Let $\varphi_{1} ,\varphi_{2} \in O(X, \tau)$.\newline

1) If $\varphi_{1} O(X)$ is closed under finite intersection and $%
\varphi_{1} O(X) \subset \varphi_{2} O(X)$, then for each $a\in X$, the
family $\Phi _a = \varphi_{1} O(X,a)$ is a filterbase, and $\Phi _a$ $%
\varphi_{1,2}$-converges to a.\newline

2) If $\varphi_{2}$ is regular w.r.t. $\varphi_{1} O(X)$ and $\varphi_{1}
O(X) \subset \varphi_{2} O(X)$, then for each $a\in X$ the family $\Phi _a
=\{\varphi_{2} (U): U\in \varphi_{1} O(X,a)\}$ is a filterbase and $%
\varphi_{1,2}$-converges to a.
\end{teo}

\begin{teo}
([20]) Let $A\subset X$ and $a\in X$.\newline

1) If there exists a filter which contains A and $\varphi_{1,2}$-accumulates
to a, then $a\in \varphi_{1,2} cl A$. \newline

2) If $\varphi_{2}$ is regular w.r.t. $\varphi_{1} O(X)$ and $a\in
\varphi_{1,2} cl A$, then there exists a filter containing A and $%
\varphi_{1,2}$-converging to a.\newline

3) If $\varphi_{2}$ is regular w.r.t. $\varphi_{1} O(X)$, then $a\in
\varphi_{1,2} cl A$ iff there exists a filter ${\cal F}$ containing A and $%
\varphi_{1,2}$-converging to a.\newline

4) If $\varphi_{2}$ is regular w.r.t. $\varphi_{1} O(X)$, then A is $%
\varphi_{1,2}$-closed iff whenever there exists a filter containing A and $%
\varphi_{1,2}$-converging to a point a in X, then $a\in A$.
\end{teo}

{\bf Proof.} 2) Let $a\in \varphi_{1,2} clA$. Then $\varphi_{2} (U) \cap A
\neq \emptyset $ for each $U \in \varphi_{1} O(X,a)$.

$\Phi =\{\varphi_{2} (U) \cap A : U \in \varphi_{1} O(X,a)\}$ is a
filterbase. The filter generated by $\Phi$ contains A and $\varphi_{1,2} $-$%
converges$ to a.

The proofs of the others are easy.

\begin{teo}
Let us define $cl^* :P(X) \rightarrow P(X)$ as $cl^* A =\{x:$ there exists a
filter ${\cal F}$ containing A and $\varphi_{1,2}$-converging to x$\}$ for $%
A\subset X$.\newline

1) If $\varphi_{2}$ is regular w.r.t. $\varphi_{1} O(X)$ then $cl^*
A=\varphi_{1,2} clA$ and $cl^*$ operator defines the following (pre)topology.

$\tau ^* =\{U\subset X : cl^*(X\setminus U) \subset X\setminus
U\}=\{U\subset X: \varphi_{1,2} cl (X\setminus U)\subset X\setminus U\}=\tau
_{\varphi_{1,2}}$. \newline

2) If $\varphi_{2}$ is regular w.r.t. $\varphi_{1} O(X)$ and $\varphi_{1}
O(X) \subset \varphi_{2} O(X)$ then $cl^*$ operator defines the following
topology.

$\tau ^* =\{U\subset X : cl^*(X\setminus U)=( X\setminus U)\}=\{U\subset X:
\varphi_{1,2} cl (X\setminus U)= X\setminus U\}=\tau _{\varphi_{1,2}}$.
\newline

3) If $\varphi_{2}$ is regular w.r.t. $\varphi_{1} O(X)$, $\varphi_{1} O(X)
\subset \varphi_{2} O(X)$, and for each $U\in \varphi_{1} O(X)$ $\varphi_{2}
(U) \in \varphi_{1,2} O(X)$, then $cl^*$ operator is a Kuratowski closure
operator and again $\tau ^* =\tau _{\varphi_{1,2}}$.
\end{teo}

\begin{ex}
Let $\varphi_{1} = int$,\quad $\varphi_{2} =cl$.

A filter ${\cal F}$ $\varphi_{1,2}$-converges to a iff ${\cal F}$ is $\theta$%
-converges to a. $\theta$-convergence defines a pretopology [8].
\end{ex}

\section{$\protect\varphi_{1,2}$-compactness}

\begin{defi}
([20]) Let $\varphi_{1} ,\varphi_{2} \in O(X, \tau)$, $X\in {\cal A} \subset
P(X)$, $A\subset X$.\newline

a) If $A\subset \cup {\cal U}$ for a subfamily ${\cal U} $ of ${\cal A}$,
then ${\cal U}$ is called an ${\cal A}$-cover of A. If an ${\cal A}$-cover $%
{\cal U}$ of A is countable (finite) then we call ${\cal U}$ as countable $%
{\cal A}$-cover (finite ${\cal A}$-cover).\newline

b) If each ${\cal A}$-cover ${\cal U}$ of A has a finite subfamily ${\cal U}
^{^{\prime}}$ such that $A\subset \cup \{\varphi_{2} (U): U\in {\cal U}
^{^{\prime}}\}$, then we call A is $({\cal A}$-$\varphi_{2} )$-compact
relative to X (shortly (${\cal A}$-$\varphi_{2}$)-compact set).\newline

c) We call an (${\cal A}$-$\imath$)-compact set relative to X as ${\cal A}$%
-compact set shortly.\newline

d) If we take ${\cal A} =\varphi_{1} O(X)$ in (b), then we call A as $%
\varphi_{1,2}$-compact relative to X (shortly $\varphi_{1,2}$-compact set).

If we take ${\cal A} =\varphi_{1,2} O(X)$ in (c) we get the definition of a $%
\varphi_{1,2} O(X)$-compact set.

If X is $\varphi_{1,2}$-compact set relative to itself, then X will be
called $\varphi_{1,2}$-compact space.

If X is $\varphi_{1,2} O(X)$-compact set relative to itself, then X will be
called $\varphi_{1,2} O(X)$-compact space.
\end{defi}

$\varphi_{1,2}$-Lindel\"of sets relative to X, $\varphi_{1,2}$-Lindel\"of
spaces and $\varphi_{1,2}$-countable compact sets relative to X, $%
\varphi_{1,2}$-countable compact spaces were defined in a similar way as in
[21], [22].

\begin{teo}
([20]) If $\varphi_{1}^{^{\prime}} O(X) \subset \varphi_{1} O(X)$, $%
\varphi_{2} \leq \varphi_{2}^{^{\prime}}$ (hence if $\varphi_{1}^{^{\prime}}
\leq \varphi_{1} , \varphi_{2} \leq \varphi_{2}^{^{\prime}}$), then each $%
\varphi_{1,2}$-compact set is $\varphi_{1,2} ^{^{\prime}}$-compact set.
\end{teo}

\begin{ex}
Let $A\subset X$. \newline

a) Let $\varphi_{1} =int$, $\varphi_{2} =cl$.\newline
A is $\varphi_{1,2}$-compact set iff A is an H-set.\newline

b) Let $\varphi_{1} =int$, $\varphi_{2} = intocl$.\newline
A is $\varphi_{1,2}$-compact set iff A is an N-set.\newline

c) Let $\varphi_{1} =cloint$, $\varphi_{2} = scl$.\newline
A is $\varphi_{1,2}$-compact set iff A is an s-set.\newline

d) Let $\varphi_{1} = cloint$, $\varphi_{2} =cl$.\newline
A is $\varphi_{1,2}$-compact set iff A is an S-set.\newline

e) Let $\varphi_{1} =int$, $\varphi_{2} =\imath$.\newline
A is $\varphi_{1,2}$-compact set iff A is compact.
\end{ex}

By using Theorem 3.2, we get that, each N-set is an H-set, each s-set is an
S-set and each S-set is an H-set.

\begin{teo}
([20]) The following are equivalent for any subset A of X.\newline

a) A is $\varphi_{1,2}$-compact set.\newline

b) Every filterbase in X which meets A, $\varphi_{1,2}$-accumulates in X to
some point in A.\newline

c) Every maximal filterbase in X which meets A, $\varphi_{1,2}$-converges in
X to some point in A. \newline

d) Every filterbase in A $\varphi_{1,2}$-accumulates in X to some point in A.%
\newline

e) Every maximal filterbase in A, $\varphi_{1,2}$-converges in X to some
point in A.\newline

f) For any family W of nonempty sets with $A\cap (\cap \{\varphi_{1,2} cl F
:F\in W\})=\emptyset$, there exists a finite subfamily $W^{^{\prime}}$ of W
such that $A\cap (\cap W^{^{\prime}})=\emptyset$. \newline

g) For any family of nonempty sets such that for each finite subfamily $%
W^{^{\prime}}$ of W we have $A\cap (\cap \{F: F\in W^{^{\prime}} \})\neq
\emptyset$, then $A\cap (\cap \{\varphi_{1,2} cl F :F\in W\})\neq \emptyset$%
\newline

h) If ${\cal F}$ is a filterbase such that $A\cap \{\varphi_{1,2} cl F :F\in
{\cal F} \}=\emptyset$ then there exists an $F\in {\cal F}$ such that $F\cap
A =\emptyset$.

If $\tilde \varphi_{2}$ is the dual of $\varphi_{2}$, then the following
statements (i) and (j) are equivalent to each one of the above statements.%
\newline

i) For any family $\Phi$ of $\varphi_{1}$-closed sets with $A\cap (\cap
\Phi) =\emptyset$, there exists a finite subfamily $\Phi ^{^{\prime}}$ of $%
\Phi$ such that $A\cap (\cap \{\tilde \varphi_{2} (F): F\in \Phi
^{^{\prime}} \})=\emptyset$. \newline

j) If $\Phi$ is a family of $\varphi_{1}$-closed sets such that for each
finite subfamily $\Phi ^{^{\prime}}$ of $\Phi$ we have $A\cap (\cap \{\tilde
\varphi_{2} (F): F\in \Phi ^{^{\prime}} \})\neq \emptyset$, then $A\cap
(\cap \Phi)\neq \emptyset$.
\end{teo}

By choosing X instead of A in the above Theorem 3.4, we get the equivalent
statements for a space $(X, \tau)$ to be $\varphi_{1,2}$-compact space.

\begin{teo}
Let ${\cal B} =\{\varphi_{2} (U): U\in \varphi_{1} O(X)\}$.\newline

a) If $\varphi_{2} (U)\in \varphi_{1} O(X)$, $\varphi_{2} (\varphi_{2}
(U))\subset \varphi_{2} (U)$ for each $U\in \varphi_{1} O(X)$ then ${\cal B}
\subset \varphi_{1,2} O(X) \cap \varphi_{1} O(X)$.\newline

b) If $\varphi_{1} O(X) \subset \varphi_{2} O(X)$ and if ${\cal B} \subset
\varphi_{1,2} O(X)$, then ${\cal B}$ is a base for the supratopology $%
\varphi_{1,2} O(X)$.\newline

c) If ($\varphi_{2} \geq \varphi_{1}$ or $\varphi_{2} \geq \imath$) and if $%
{\cal B} \subset \varphi_{1,2} O(X)$, then ${\cal B}$ is a base for the
supratopology $\varphi_{1,2} O(X)$ [17].\newline

d) If ($\varphi_{2} \geq \varphi_{1}$ or $\varphi_{2} \geq \imath$) and, $%
\varphi_{2} (U) \in \varphi_{1} O(X)$, $\varphi_{2}(\varphi_{2} (U))\subset
\varphi_{2} (U)$ for each $U\in \varphi_{1} O(X)$, then ${\cal B}$ is a base
for the supratopology $\varphi_{1,2} O(X)$ [18].
\end{teo}

{\bf Proof.} a) Let $U\in \varphi_{1} O(X)$ and $x\in \varphi_{2} (U)$. $%
x\in \varphi_{2} (U) \in \varphi_{1} O(X)$ and $\varphi_{2} (\varphi_{2}
(U))\subset \varphi_{2}(U)$. So $x\in \varphi_{1,2} int \varphi_{2} (U)$. We
have $\varphi_{2} (U) \subset \varphi_{1,2} int \varphi_{2} (U)$. Hence $%
\varphi_{2} (U) \in \varphi_{1,2} O(X)$ for each $U\in \varphi_{1} O(X)$ and
${\cal B} \subset \varphi_{1,2} O(X)\cap \varphi_{1} O(X)$.\newline

b) Let $A \in \varphi_{1,2} O(X)$ and $x\in A$. There exists a $U\in
\varphi_{1} O(X,x)$ such that $\varphi_{2} (U) \subset A$. We have $x\in
U\subset \varphi_{2} (U) \subset A$, $\varphi_{2} (U) \in \varphi_{1,2} O(X)$
and $\varphi_{2} (U) \in {\cal B}$.

Proofs of (c) and (d) are clear from Lemma 2.1 and (a), (b).

\begin{teo}
Let ${\cal B} =\{\varphi_{2} (U): U\in \varphi_{1} O(X)\}$. If ($\varphi_{2}
\geq \varphi_{1}$ or $\varphi_{2} \geq \imath$), and for each $U\in
\varphi_{1} O(X)$ we have $\varphi_{2} (U) \in \varphi_{1} O(X)$, $%
\varphi_{2} (\varphi_{2} (U))\subset \varphi_{2} (U)$ then the following are
equivalent for any subset A of X.\newline

a) A is $\varphi_{1,2}$-compact set.\newline

b) A is ${\cal B}$-compact set.\newline

c) A is $\varphi_{1,2} O(X)$-compact set.\newline

d) If W is any subfamily of $\{X\setminus \varphi_{2} (U): U\in \varphi_{1}
O(X)\}$ such that for each finite subfamily $W ^{^{\prime}}$ of W we have $%
A\cap (\cap W ^{^{\prime}})\neq \emptyset$ then $A\cap (\cap W)\neq \emptyset
$.\newline

e) If W is any subfamily of $\{X\setminus \varphi_{2} (U) : U\in \varphi_{1}
O(X) \}$ with $A\cap (\cap W)= \emptyset$ then there exists a finite
subfamily $W ^{^{\prime}}$ of W such that $A\cap (\cap W ^{^{\prime}})=
\emptyset$.\newline

f) If W is any subfamily of $\{X\setminus U: U\in \varphi_{1,2} O(X)\}$ such
that for each finite subfamily $W ^{^{\prime}}$ of W we have $A\cap (\cap W
^{^{\prime}})\neq \emptyset$ then $A\cap (\cap W)\neq \emptyset$.\newline

g) If W is any subfamily of $\{X\setminus U : U\in \varphi_{1,2} O(X)\}$
such that $A \cap (\cap W)= \emptyset$, then there exists a finite subfamily
$W ^{^{\prime}}$ of W such that $A \cap (\cap W ^{^{\prime}})=\emptyset$.
\end{teo}

{\bf Proof.} Let us see that $a\Leftrightarrow c$ is true.

In the case $A= \emptyset$ proofs are clear.

Let A be a nonempty $\varphi_{1,2}$-compact set and ${\cal U}$ a $%
\varphi_{1,2} O(X)$-cover of A. For each $x\in A$, there exists a $U_x \in
{\cal U}$ s.t. $x\in U_x$. There exists a $\varphi_{1}$-open set $V_x$
containing x s.t. $V_x \subset \varphi_{2} (V_x) \subset U_x$ .Since A is $%
\varphi_{1,2}$-compact set, there exists a finite subset $\{x_1,...,x_n\}$
of A s.t. $A\subset \cup^{n}_{i=1}\varphi_{2} (V_{x_i})\subset
\cup^{n}_{i=1} U_{x_i} $. Hence A is $\varphi_{1,2} O(X)$-compact set.

Let A be $\varphi_{1,2} O(X)$-compact set and ${\cal U}$ a $\varphi_{1} O(X)$%
-cover of A. We have $A\subset \cup \{\varphi_{2} (U) : U\in {\cal U} \}$.
Since for each $U\in {\cal U}$, $\varphi_{2} (U) \in {\cal B} \subset
\varphi_{1,2} O(X)$, there exists a finite subfamily $\{U_1,...,U_n\}$ of $%
{\cal U}$ s.t. $A\subset \cup ^{n}_{i=1} \varphi_{2} (U_i)$. Hence A is $%
\varphi_{1,2}$-compact set.

Under the given conditions, since ${\cal B}$ is a base of the supratopology $%
\varphi_{1,2} O(X)$, $b\Leftrightarrow c$ is clear. Now the other proofs are
easy.

Under the hypothesis of Theorem 3.6 by joining Theorems 3.4 and 3.6 we get
the equivalent statements for a set to be $\varphi_{1,2}$-compact set.

\begin{teo}
Under the hypothesis of Theorem 3.6, the following are equivalent.

a) X is $\varphi_{1,2}$-compact space.\newline

b) X is ${\cal B}$-compact space.\newline

c) X is $\varphi_{1,2} O(X)$-compact space.\newline

d) For each $U\in \varphi_{1} O(X)$, $X\setminus \varphi_{2} (U)$ is $%
\varphi_{1,2}$-compact set.\newline

e) For each $U\in \varphi_{1} O(X)$, $X\setminus \varphi_{2} (U)$ is $%
\varphi_{1,2} O(X)$-compact set.\newline

f) For each $U\in \varphi_{1} O(X)$, $X\setminus \varphi_{2} (U)$ is ${\cal B%
}$-compact set.\newline

g) Each $\varphi_{1,2}$-closed set is $\varphi_{1,2}$-compact set.\newline

h) Each $\varphi_{1,2}$-closed set is $\varphi_{1,2} O(X)$-compact set.%
\newline

i) Each $\varphi_{1,2}$-closed set is ${\cal B}$-compact set.
\end{teo}

Now, by using the Theorems 3.4, 3.6 and 3.7, we get the equivalent forms for
a space $(X, \tau)$ to be $\varphi_{1,2}$-compact space under the hypothesis
of Theorem 3.6.

\begin{ex}
Let $A\subset X$.\newline

a) Let $\varphi_{1} =int$, $\varphi_{2} =intocl$ as in Example 3.3 (b).

$\varphi_{1} O(X) =\tau$. $\varphi_{2} \geq \varphi_{1}$. For each $U\in
\varphi_{1} O(X) =\tau$, we have $\varphi_{2} (U) =U^{\frac {o}{}}\in
\varphi_{1} O(X)$ and $\varphi_{2} (\varphi_{2} (U)) =(U^{\frac {o}{}})^{%
\frac {o}{}}=U^{\frac {o}{}}=\varphi_{2} (U)$. ${\cal B} =\{\varphi_{2} (U)
:U\in \varphi_{1} O(X)\}= \{U^{\frac {o}{}}: U\in \tau \} =RO(X)$. $%
\varphi_{1,2} O(X) = \tau _s$. We know that RO(X) is a base for $\tau _s$. $%
\{X\setminus \varphi_{2} (U) : U\in \varphi_{1} O(X)\}=RC(X)$. $%
\varphi_{1,2} C(X) =$ the family of $\tau _s$-closed sets = the family of $%
\delta$-closed sets. $\tilde \varphi_{2} = cloint$ is the dual of $%
\varphi_{2}$.

A is $\varphi_{1,2}$-compact set iff it is N-set iff A is compact in $(X,
\tau_s)$ iff A is RO(X)-compact set.

$(X, \tau)$ is $\varphi_{1,2}$-compact iff ($X, \tau _s$) is compact. These
are very well known results.\newline

b) Let $\varphi_{1} =cloint$, $\varphi_{2} =scl$ as in Example 3.3 (c).%
\newline
$\varphi_{2} \geq \imath$. $\varphi_{1} O(X) = SO(X)$, $\varphi_{1,2} O(X) =
S\theta O(X)$, $\varphi_{1,2} C(X) = S\theta C(X).$ For each $U\in
\varphi_{1} O(X)$ we have $\varphi_{2} (U) = scl U \in SR(X) \subset SO(X)$,
and $\varphi_{2} (\varphi_{2} (U))= scl (scl U) = sclU =\varphi_{2} (U)$. $%
{\cal B} =\{\varphi_{2} (U) : U\in \varphi_{1} O(X)\}= SR(X)$ is a base for
the supratopology $S\theta O(X)$. $\{X\setminus \varphi_{2} (U): U\in
\varphi_{1} O(X)\}=SR(X)$. $\tilde \varphi_{2} =$semi-int is the dual of $%
\varphi_{2}$.

A is $\varphi_{1,2}$-compact set iff it is (SO(X)-scl)-compact set iff it is
$S\theta O(X)$-compact set iff it is SR(X)-compact set.
\end{ex}

Now, by using the Theorems 3.4, 3.6 we can write the equivalent statements
for a set to be ($\tau$-intocl)-compact set, or to be (SO(X)-scl)-compact
set, and by using the Theorems 3.4, 3.6, 3.7 we can write the equivalent
statements for a space to be ($\tau$-intocl)-compact space or to be
(SO(X)-scl)-compact space.

Some equalities related to closure types can be obtained by using operations
and some of them were given in [18]. For example for any open set T in a
topological space $(X, \tau)$, we have $\overline T =\theta clT = \tau _s$-$%
clT$.

\begin{teo}
Let $A\subset X$. If $\varphi_{2} (U) =\varphi _3 (U)$ for each $U\in
\varphi_{1} O(X)$, then $\varphi_{1,2} O(X) = \varphi _{1,3} O(X)$, and A is
$\varphi_{1,2}$-compact ($\varphi _{1,2} O(X)$-compact ) set iff it is
\newline
$\varphi_{1,3} $-compact ($\varphi _{1,3} O(X)$-compact) set.
\end{teo}

\begin{ex}
Let $\varphi_{1} =int$, $\varphi_{2} =scl$, $\varphi_3 =intocl$.

For $U\in \varphi_{1} O(X)=\tau$, $\varphi_{2} (U) =sclU =U\cup U^{\frac {o}{%
}}=U^{\frac {o}{}} =\varphi _3 (U)$, (sclU = $U\cup U^{\frac {o}{}}$, [3]).
A is $\tau_s$-compact set iff A is ($\tau$-intocl)-compact set iff A is
RO(X)-compact set iff A is ($\tau$-scl)-compact set.
\end{ex}

\begin{teo}
If $\varphi_{1}$ is monotonous and, for each pair $U, V \in \varphi_{1} O(X)$%
, $\varphi_{2} (U\cup V)= \varphi_{2} (U) \cup \varphi_{2} (V)$, then the
following are equivalent.\newline

a) A is $\varphi_{1,2}$-compact set.\newline

b) Each filterbase ${\cal F} _b$ which is a subfamily of $\{X\setminus
\varphi_{2} (U): U\in \varphi_{1} O(X)\}$ and meets A, $\varphi_{1,2}$%
-accumulates to some point $a\in A$.
\end{teo}

{\bf Proof.} (a$\Rightarrow$b). It is clear from Theorem 3.4.

(b$\Rightarrow$a). Let A be not $\varphi_{1,2}$-compact set under the
assumption of (b). There is a subfamily ${\cal U} =\{U_i : i\in I\}$ of $%
\varphi_{1} O(X)$ s.t. $A \subset \cup {\cal U}$ but for each finite subset
J of I $A\not\subset \cup _{j\in J}\varphi_{2} (U_j)= \varphi_{2} (\cup
_{j\in J}U_j)$. For each finite subset J of I we have $A \cap (X\setminus
\varphi_{2} (\cup _{j\in J}U_j))\neq \emptyset$. Since $\varphi_{1}$ is
monotonous, $\varphi_{1} O(X)$ is a supratopology and hence ${\cal F} _b
=\{X\setminus \varphi_{2} (\cup _{j\in J}U_j): J\subset I, J finite\}$ is a
filterbase s.t. ${\cal F} _b \subset \{X \setminus \varphi_{2} (U) : U\in
\varphi_{1} O(X)\}$ and ${\cal F} _b$ meets A. There exists a point a in A
s.t. ${\cal F} _b$ $\varphi_{1,2}$-accumulates to a. There exists a $U_{i_a}
\in {\cal U}$ s.t. $a\in U_{i_a}$. $X\setminus \varphi_{2} (U_{i_a}) \in
{\cal F} _b$, $\varphi_{2} (U_{i_a}) \cap (X\setminus \varphi_{2}
(U_{i_a}))\neq \emptyset$. This contradiction completes the proof.

\begin{ex}
a) Let $A\subset X$, $\varphi_{1} = cloint$, $\varphi_{2} = cl$ as in
Example 3.3 (d).

$\varphi_{1} O(X) =SO(X)$, $\varphi_{2} \geq \varphi_{1}$. $\varphi_{1}$ is
monotonous. For each $U\in \varphi_{1} O(X) =SO(X)$, $\varphi_{2} (U) =
\overline U\in SO(X)$ and $\varphi_{2} (\varphi_{2}(U))= \varphi_{2} (U)$.
If $U,V \in \varphi_{1} O(X)$, then $\varphi_{2} (U\cup V)= \overline {U\cup
V} =\overline U \cup \overline V= \varphi_{2} (U) \cup \varphi_{2} (V)$. $%
{\cal B} =\{\varphi_{2} (U) : U\in \varphi_{1} O(X)\}= RC(X)$. $%
\varphi_{1,2} O(X)=\theta SO(X)$.

A is $\varphi_{1,2}$-compact set iff it is ${\cal B}$-compact set iff it is $%
\varphi_{1,2} O(X)$-compact set. X is $\varphi_{1,2}$-compact and Hausdorff
iff X is S-closed space.

$\tilde \varphi_{2} = int$ is the dual of $\varphi_{2}$. $\{X\setminus
\varphi_{2} (U) : U\in \varphi_{1} O(X)\} =RO(X) = \{\tilde \varphi_{2}
(X\setminus U): U\in \varphi_{1} O(X)\}= \{int K: K\in SC(X)\}$.

b) Let $\varphi_{1} = int$, $\varphi_{2} = cl$ as in Example 3.3 (a).

$\varphi_{1}$ is monotonous, $\varphi_{2} (U\cup V) = \varphi_{2} (U) \cup
\varphi_{2} (V)$ for $U, V \in \varphi_{1} O(X) =\tau$. $\{X\setminus
\varphi_{2} (U) : U\in \varphi_{1} O(X)\} =RO(X)$. For a filterbase ${\cal F}
_b \subset RO(X)$, we have $\cap \{\varphi_{1,2} cl F : F \in {\cal F}
_b\}=\cap \{\theta cl F : F\in {\cal F} _b\}=\cap \{\overline F : F\in {\cal %
F} _b\}= \cap \{\tau _s $-$clF: F\in {\cal F} _b\}$. X is $\varphi_{1,2}$%
-compct and Hausdorff iff X is H-closed space.
\end{ex}

Now we can get many known results that some of them can be seen from
[4,5,7,8] and many unknown results by special choices of operations. And we
will have many results by using the papers [18,19,20].

\begin{center}
{\bf REFERENCES}
\end{center}

\begin{enumerate}
\item[1.] M.E.Abd El-Monsef, F.M.Zeyada and A.S.Mashhour, Operations on
topologies and their applications on some types of coverings. Ann. Soc. Sci.
Bruxelles, 97, No. 4 (1983), 155-172.

\item[2.] M.E.Abd El-Monsef and E.F.Lashien, Local discrete extensions of
supratopologies. Tamkang J. Math. 21, No.1(1990), 1-6.

\item[3.] D. Andrijevi\'{c}, On the topology generated by pre-open sets,
Mat. Vesnik, 39 (1987), 367-376.

\item[4.] R. F. Dickman, JR and Jack R. Porter, $\Theta$-perfect and
absolutely closed functions, Illinois J. Math. 21 (1977),42-60.

\item[5.] R. F. Dickman, JR and R. L. Krystock, S-sets and s-perfect
mappings, Proc. Amer. Math. Soc. 80, 4 (1980), 687-692.

\item[6.] J. Dugundji, Topology, (1966).

\item[7.] L. L. Herrington and P. E. Long, Chracterizations of H-closed
spaces, Proc. Amer. Math. Soc. 48, 2 (1975), 469-475.

\item[8.] R. A. Herrmann, RC-Convergence, Proc. Amer. Math. Soc. 75, 2
(1979), 311-317.

\item[9.] D. S. Jankovi\'{c}, On functions with $\alpha$-closed graphs,
Glasnik Matemat{\i}\v{c}k{\i}, 18, No. 38(1983), 141-148.

\item[10.] A. Kandil, E. E. Kerre and A. A. Nouh, Operations and mappings on
fuzzy topological spaces, Ann. Soc. Sci. Bruxelles  105, No. 4. (1991),
165-188.

\item[11.] S. Kasahara, Operation-compact spaces, Math. Japonica 24, No. 1
(1979), 97-105.

\item[12.] E. E. Kerre, A. A. Nouh and A. Kandil, Operations on the class of
fuzzy sets on a universe endowed with a fuzzy topology, J. Math. Anal. Appl.
180 (1993), 325-341.

\item[13.] J. K. Kohli, A class of mappings containing all continuous and
all semiconnected mappings, Proc. Amer. Math. Soc. 72, No. 1 (1978), 175-181.

\item[14.] J. K. Kohli, A framework including the theories of continuous and
certain non-continuous functions, Note di Mat. 10, No. 1 (1990), 37-45.

\item[15.] J. K. Kohli, Change of topology, characterizations and product
theorems for semilocally P-spaces, Houston J. Math.  17, No. 3 (1991),
335-349.

\item[16.] M. N. Mukherjee and G. Sengupta, On $\pi$-closedness: A unified
theory, Anal. Sti. Univ. "Al. I. Cuza" Iasi, 39, s.l.a, Mat. No. 3 (1993),
1-10.

\item[17.] G. Sengupta and M. N. Mukherjee, On $\pi$-closedness a unified
theory II, Boll. U. M. I. 7, 8-B (1994), 999-1013.

\item[18.] T. H. Yalva\c{c}, A unified theory for continuities, submitted.

\item[19.] T. H. Yalva\c{c}, A unified theory on openness and closedness of
functions, submitted.

\item[20.] T. H. Yalva\c{c}, A unified approach to compactness and filters,
appear.

\item[21.] T. H. Yalva$\d{c}$, On some unifications, The First Turkish
International Conference on Topology and Its Applications, August (2000),
submitted.

\item[22.] T. H. Yalva\c{c}, Unification of some concepts similar to
Lindel\"of property , submitted.
\end{enumerate}

\end{document}